\documentclass[a4paper,12pt]{article}
\usepackage[T2A]{fontenc}
\usepackage{amsfonts}
\usepackage{amsmath,amssymb,amsthm}
\usepackage{wrapfig}

\newcommand{\vol}{\mathop{\rm vol}\nolimits}
\newcommand{\ar}{\mathop{\rm area}\nolimits}

\title{On the Optimality of Functionals \\ over Triangulations of Delaunay Sets}

\author{Nikolay P. Dolbilin,
        Herbert Edelsbrunner and
        Oleg R. Musin\thanks{This research is supported by the Russian
          government project 11.G34.31.0053,
          RFBR grant 11-01-00735, and DMS 1101688.}}

\begin{document}
\date{}
\maketitle

In this short paper, we consider the functional density on sets
of uniformly bounded triangulations with fixed sets of vertices.
We prove that if a functional attains its minimum on
the Delaunay triangulation, for every finite set in the plane,
then for infinite sets the density of this functional attains
its minimum also on the Delaunay triangulations.

A {\it Delaunay set} in $\Bbb{R}^d$ is a set of points, $X$,
for which there are positive $r$ and $R$ such that
every open $d$-ball of radius $r$ contains at most one point
and every closed $d$-ball of radius $R$ contains at least one point of $X$.
In this paper, we consider Delaunay sets in general position,
i.e.\ no $d+2$ points in $X$ lie on a common $(d-1)$-sphere. 
By a {\it triangulation} of $X$ we mean a simplicial complex whose
vertex set is $X$.
For finite sets, the simplices decompose the convex hull of the set,
while for Delaunay sets $X$, the simplices decompose $\Bbb{R}^d$.
We say a triangulation $T$ is {\it uniformly bounded} if there exists
a positive number $q = q(T)$ that is larger than or equal
to the circumradii of all $d$-simplices in the triangulation:
${\cal R}(S) \le q$ for all $d$-simplices $S$ of $T$.
We denote the family of all uniformly bounded triangulations
of $X$ by $\Theta(X)$.

Delaunay sets have been introduced by Boris Delaunay in 1924,
who called them {\it $(r,R)$-systems}.
He proved that for any Delaunay set $X$, there exists a unique
Delaunay tessellation $DT(X)$; see e.g.\ \cite{Del}.
If $X$ is in general position,
then $DT(X)$ is a triangulation of $X$ in the sense defined above.
Since the circumradius of any simplex is at most $R$,
the Delaunay triangulation is uniformly bounded with $q = R$,
i.e.\ $DT(X) \in \Theta (X)$.
We note that every Delaunay set also has triangulations that are not
uniformly bounded, and it is not difficult to contruct them.

We want to remind the reader of a related open problem
about Delaunay sets:
{\it is it true that for every planar Delaunay set, $X$,
and for every positive number, $C$,
there exists a triangle $\Delta$ that contains none of the points in $X$
and whose area is greater than $C$}?
While we heard of this question from Michael Boshernitzan, it is sometimes
referred to as {\it Danzer's problem}.

Let $F$ be a functional defined on $d$-simplices $S$.
For instance, $F(S)$ may be the sum of squares of edge lengths
multiplied by the volume of $S$.
We only consider functionals that are continuous in parameters
describing the simplices, for example the lengths of their edges.
Let $Y$ be a finite set in ${\Bbb{R}}^d$ and $T$ any triangulation of $Y$.
Then $F$ can be defined on $T$ as
\begin{eqnarray*}
  F(T)  &=&  \sum_{S \in T} F(S).
\end{eqnarray*}
It is clear, that this definition cannot be used for infinite sets.
We therefore define the {\it (lower) density} of $F$ for a uniformly bounded
triangulation $T$ of a Delaunay set $X$ as
\begin{eqnarray*}
  {\bar F}(T)  &=&  \mbox{\rm lower}\lim_{\alpha \rightarrow \infty}
                    {\frac{1}{\vol(B_\alpha)} \sum_{S \subset B_\alpha} F(S)},
\end{eqnarray*}
where $B_\alpha$ denotes the closed ball of radius $\alpha$ with center at the
origin of $\Bbb{R}^d$.
We call ${\bar F}(T)$ the {\it (lower) density} of $F$ on $T$.
For the rest of the paper, we limit ourselves to dimension $d = 2$.

\medskip
\noindent\textbf{Theorem.} \textit{Let $F$ be a continuous functional
  that attains its minimum for every finite set $Y \subset {\Bbb R}^2$
  for the Delaunay triangulation of $Y$.
  Then the density $\bar F$ on $\Theta (X)$, where $X \subset {\Bbb R}^2$
  is a Delaunay set, attains its minimum for the Delaunay triangulation of $X$.}

\begin{proof}
 Let $T \in \Theta (X)$ be a triangulation with parameter $q$,
 and consider the simplicial complex $T_\alpha \subset T$
 that consists of all triangles, edges, and vertices of $T$
 contained in $B_{\alpha}$.
 We consider the convex hull of the vertices of $T_\alpha$,
 which we denote as $C_\alpha$.
 The difference between $C_\alpha$ and the union of triangles in $T_\alpha$
 consists of polygons,
 and since any polygon can be triangulated without adding vertices,
 $T_\alpha$ can be extended to a triangulation $T_\alpha'$
 of the same set of vertices.
 Write $K_\alpha$ for the number of triangles in $T_\alpha$.
 Since the circumradius of each triangle is bounded from above,
 and the lengths of its edges are bounded from below,
 the area of each triangle is at least some constant.
 It follows that $K_\alpha$ is at most some contant times $\alpha^2$.
 The circle bounding $B_\alpha$ intersects at most some constant times
 $\sqrt{K_\alpha}$ of the triangles in $T$, which implies
 that $T_\alpha'$ has at most some constant times $\sqrt{K_\alpha}$
 triangles in addition to those in $T_\alpha$.
 Using the continuity of the functional, 
 it follows that the limit of the ratio $F (T_\alpha)$ over $F (T_\alpha')$,
 for $\alpha$ going to infinity, is equal to $1$.
 By the assumption, $F (T_\alpha')$ is no less than the value of
 $F$ on the Delaunay triangulation of the same set of vertices,
 which completes the proof.
\end{proof}

We remark that there are non-convex polytopes in dimension $d > 2$
that cannot be triangulated without adding new vertices.
They constitute the main difficulty in extending the Theorem
to general dimension.
The theorem and the results stated in the papers
\cite{Del,Lam,Mus1,Mus2,Raj} yield

\medskip

\noindent\textbf{Corollary}. 
{\it Let $\Delta$ be a triangle with barycenter $b$, circumcenter $c$,
and edges of lengths $a_1, a_2, a_3$.
     Let us consider the following functionals:\\
 (1) $F_1(\Delta) = {\cal R}^a(\Delta)$,
     where ${\cal R}(\Delta)$ is the circumradius and $a>0$;\\
 (2) $F_2(\Delta) = (a_1^2+a_2^2+a_3^2)/{\ar(\Delta)}$,
     where $\ar(\Delta)$ is the area of $\Delta$;\\
 (3) $F_3(\Delta) = -\rho(\Delta)$,
     where $\rho (\Delta)$ is the inradius of $\Delta$; \\
 (4) $F_4(\Delta) = (a_1^2+a_2^2+a_3^2) \ar(\Delta)$;\\
 (5) $F_5(\Delta) = {\cal R}(\Delta)^a \ar(\Delta)$, where $a \ge 1$;\\
 (6) $F_6(\Delta) = \|b(\Delta)-c(\Delta )\|^2\ar(\Delta)$.
Then the densities ${\bar F}_i$, for $i = 1, 2, \ldots, 6$,
achieve their minima on the Delaunay triangulations of Delaunay sets
in the plane.}

\medskip

For finite sets, the optimality of the functionals $F_1$ and $F_2$
has been shown in \cite{Mus1}, of $F_3$ it has been shown in \cite{Lam},
of $F_4$ it has been shown \cite{Raj},
and of $F_5$ and $F_6$ it has been shown in \cite{Mus2}.

\end{document}